\documentclass[final]{siamltex}

\usepackage{epsfig}
\usepackage{amsmath}
\usepackage{amssymb}
\usepackage{subfigure}

\def\Cst{{\mathbb C}}
\def\Nst{{\mathbb N}}

\def\Zst{{\mathbb Z}}
\def\lsp{{\boldsymbol\ell}}
\def\ltsp{{\lsp^2}}
\def\ltZ{{\lsp^2(\Zst})}
\def\eps{\varepsilon}
\def\Csp{{\boldsymbol C}}

\def\cond{\mbox{cond\/}}
\def\conj#1{{\overline#1}}
\def\ntoinf{{n \rightarrow \infty }}

\def\Csph{{\Csp_{[-\frac{1}{2},\frac{1}{2})}}}

\newcommand{\fmin}{f_{\min}}
\newcommand{\fmax}{f_{\max}}

\newcommand{\flm}{f(\frac{l}{2m-1})}
\newcommand{\fmlm}{f_m(\frac{l}{2m-1})}
\newcommand{\finvx}{\frac{1}{f(x)}}
\newcommand{\finvlm}{\frac{1}{f(\frac{l}{2m-1})}}
\newcommand{\finvpm}{\frac{1}{f(\frac{p}{2m-1})}}
\newcommand{\fminvlm}{\frac{1}{f_m(\frac{l}{2m-1})}}

\newcommand{\explm}{e^{2\pi i kl/(2m-1)}}
\newcommand{\expkx}{e^{2\pi i kx}}
\newcommand{\gam}{\gamma}
\newcommand{\kap}{\kappa}
\newcommand{\skap}{\sqrt{\kappa}}
\newcommand{\gami}{\gamma_1}
\newcommand{\gamii}{\gamma_2}
\newcommand{\lam}{\lambda}
\newcommand{\lamkl}{\lambda^{|k-l|}}
\newcommand{\ind}{{\cal I}}

\newcommand{\polmat}{{\cal Q}_{s}}

\newcommand{\expmata}{{\cal E}_{\gamma}}
\newcommand{\expmatai}{{\cal E}_{\gamma_1}}
\newcommand{\expmatab}{{\cal E}_{\gamma,\lambda}}
\newcommand{\Ln}{L_n}
\newcommand{\Tn}{T_n}
\newcommand{\An}{A_n}
\newcommand{\Lni}{\Ln^{-1}}

\newcommand{\Ani}{\An^{-1}}
\newcommand{\Ankl}{(\An)_{kl}}
\newcommand{\Anikl}{(\Ani)_{kl}}

\newcommand{\LI}{L^{-1}}
\newcommand{\Li}{L^{-1}}

\newcommand{\Cn}{C_n}

\newcommand{\Bn}{B_n}
\newcommand{\Bna}{\Bn^{\ast}}

\newcommand{\Bnakl}{(\Bna)_{kl}}
\newcommand{\PN}{P_n}
\newcommand{\En}{E_n}
\newcommand{\EnN}{E_n^{(N)}}
\newcommand{\Sn}{S_n}
\newcommand{\Sni}{\Sn^{-1}}
\newcommand{\In}{I_n}
\newcommand{\Cnn}{C_{2n}}
\newcommand{\Cnni}{\Cnn^{-1}}
\newcommand{\egamkl}{e^{-\gam |k-l|}}

\newcommand{\egamiikl}{e^{-\gam_2 |k-l|}}

\newcommand{\egamiijl}{e^{-\gam_2 |j-l|}}
\newcommand{\egamnkl}{e^{-\gam (n- |k-l|)}}

\newcommand{\egamnkj}{e^{-\gam (n- |k-j|)}}
\newcommand{\xn}{x^{(n)}}
\newcommand{\yn}{y^{(n)}}
\newcommand{\ynn}{\tilde{y}^{(n)}}
\newcommand{\zn}{z^{(n)}}
\newcommand{\znn}{\tilde{z}^{(n)}}
\newcommand{\AI}{A^{-1}}
\newcommand{\Ai}{A^{-1}}

\newcommand{\errn}{u^{(n)}}
\newtheorem{remark}[theorem]{Remark}

\def\adots{\mathinner{\mkern1mu\raise1pt\vbox{\kern7pt\hbox{.}}\mkern2mu
   \raise4pt\hbox{.}\mkern2mu\raise7pt\hbox{.}\mkern1mu}}

\begin{document}

\title{\bf Four short stories about Toeplitz matrix calculations}
\author{Thomas Strohmer\thanks{
        Department of Mathematics, University of California, Davis, 
        CA 95616-8633, USA;  \mbox{E-mail: strohmer@math.ucdavis.edu.}
        This work was partially supported by NSF grant 9973373.}}
\date{}
\maketitle

\begin{abstract}
The stories told in this paper are dealing with the 
solution of finite, infinite, and biinfinite Toeplitz-type systems.
A crucial role plays the off-diagonal decay behavior of Toeplitz matrices
and their inverses. Classical results of Gelfand et al.\ on commutative
Banach algebras yield a general characterization of this decay behavior.
We then derive estimates for the approximate solution of (bi)infinite Toeplitz 
systems by the finite section method, showing that the approximation rate
depends only on the decay of the entries of the Toeplitz matrix and its
condition number. Furthermore, we give error estimates
for the solution of doubly infinite convolution systems by finite
circulant systems. Finally, some quantitative results on the construction of
preconditioners via circulant embedding are derived, which allow
to provide a theoretical explanation for numerical observations made by
some researchers in connection with deconvolution problems.
\end{abstract}


\noindent

\pagestyle{myheadings}
\thispagestyle{plain}
\markboth{Thomas Strohmer}{Four short stories about Toeplitz matrix calculations}

\section{Introduction}
\label{intro}

Toeplitz-type equations arise in many applications in mathematics,
signal processing, communications engineering, and statistics.
The excellent surveys~\cite{CN96,KS95} describe a number of
applications and contain a vast list of references.
The stories told in this paper are dealing with the (approximate)
solution of biinfinite, infinite, and finite hermitian positive definite 
Toeplitz-type systems. 
We pay special attention to Toeplitz-type systems with
certain decay properties in the sense that the entries of the matrix enjoy 
a certain decay rate off the diagonal. In many theoretical and practical
problems this decay is of exponential
or polynomial type. Toeplitz equations arising in image deblurring are
one example (since often the point spread function has exponential decay 
- or even stronger -  compact support)~\cite{NPT96}.
Kernels of integral equations also frequently show fast decay, leading to
Toeplitz systems inheriting this property (see e.g.~\cite{GHK94}).
Other examples include Weyl-Heisenberg frames with exponentially or
polynomially decaying window functions~\cite{Str98a} (yielding biinfinite 
block-Toeplitz systems with the same behavior when computing the so-called
dual window), as well as channel estimation problems in 
digital communications~\cite{Pro95}.

Let $\Csph$ be the set of all $1$-periodic, 
continuous, real-valued functions defined on $[-\frac{1}{2},\frac{1}{2})$.
For all $f \in \Csph$, let 
\begin{equation}
a_k=\int \limits_{-\frac{1}{2}}^{\frac{1}{2}} 
f(\omega) e^{2\pi i \omega k} \,d\omega, \qquad k=0,\pm 1, \pm 2,\dots, 
\notag
\end{equation}
be the Fourier coefficients of $f$. Since $f$ is real-valued, we have
$a_k = \conj{a}_{-k}$. 

A {\em Laurent operator} or {\em multiplication operator} associated with its
{\em defining function} $f$ can be represented
by the doubly infinite hermitian matrix $L=[L_{kl}]_{k,l=-\infty}^{\infty}$ 
with entries $L_{kl}=a_{k-l}$ for $k,l \in \Zst$. 
For all $n \ge 1$ let $\Ln = [(\Ln)_{kl}]_{k,l=-n+1}^{n-1}$ be the 
Toeplitz matrix of size $(2n-1) \times (2n-1)$ with entries 
$(\Ln)_{kl}=a_{k-l}$ for $k,l=-n+1,\dots,n-1$. $\Ln$ is a finite section
of the biinfinite Toeplitz matrix $L$.

A {\em Toeplitz operator} with {\em symbol} $f$ can be represented by 
the singly infinite matrix $T=[T_{kl}]_{k,l=0}^{\infty}$ with
$T_{kl}=a_{k-l}$ for $k,l =0,1,\dots$. 
In this case we define $\Tn= [(\Tn)_{kl}]_{k,l=0}^{n-1}$ as 
the $n \times n$ matrix with entries $(\Tn)_{kl}=a_{k-l}$ for 
$k,l=0,\dots,n-1$. Of course $\Ln = T_{2n-1}$, but in what follows
it will sometimes be convenient to use the notations $\Ln$ and $\Tn$.

As mentioned earlier, a crucial role throughout the paper plays the decay 
behavior of Toeplitz matrices and their inverses. 
Classical results of Gelfand et al.\ lead to a general characterization
of this decay behavior for biinfinite Toeplitz matrices, see 
section~\ref{s:decay}. Section~\ref{s:infinite} 
is concerned with the approximate solution of (bi)infinite Toeplitz
systems using the finite section method. Explicit error estimates are
derived, showing that the approximation rate depends only
on the condition number of the matrix and its decay properties.
In section~\ref{s:deconv} we analyze the approximate solution of
convolution equations via circulant matrices.  Finally, in 
section~\ref{s:precond}, we derive some quantitative results for
preconditioning of Toeplitz matrices by circulant embedding. Among others, we
provide a theoretical explanation of numerical observations made
by Nagy et al.\ in connection with (non)banded Toeplitz systems.

\section{On the decay of inverses of Toeplitz-type matrices} \label{s:decay}

It is helpful to review a few results on the decay of inverses of certain 
matrices. In what follows, if not otherwise mentioned, the $2$-norm of a 
matrix or a vector will be denoted by $\| . \|$ without subscript. 

The following theorem about the decay of the inverse of a band matrix
is due to Demko, Moss, and Smith~\cite{DMS84}.

\begin{theorem}
\label{th:demko}
Let $A$ be a matrix acting on $\ltsp({\cal I})$, where 
${\cal I}=\{0,1,\dots,N-1\}, \Zst$, or $\Nst$, and assume $A$ to be 
hermitian positive definite and $s$-banded (i.e., $A_{kl}=0$ if $|k-l|>s$).
Set $\kap=\|A\| \|A^{-1}\|$, $q = \frac{\skap-1}{\skap+1}$ and 
$\lam = q^{\frac{1}{s}}$. Then
\begin{equation}
|A^{-1}_{k,l}| \le c \lamkl \,,
\notag
\end{equation}
where
\begin{equation}
c = \|\AI\| \max\{1,\frac{(1+\skap)^2}{2\kap} \}\,.
\notag 
\end{equation}
\end{theorem}

One notes that the inverse of a banded matrix is in general not banded,
the type of decay changes when we switch from $A$ to $\AI$
(although exponential decay is ``almost as good'' as bandedness).
This observation suggests to look at other classes of matrices,
for which the type of decay is preserved under inversion.
This leads naturally to the following 
\begin{definition}
Let $A=[A_{k,l}]_{k,l \in \ind}$ be a matrix, where the index set is
$\ind=\Zst,\Nst$ or $\{0,\dots,N-1\}$.\\
(i) $A$ belongs to the space $\expmatab$ if the coefficients
$A_{kl}$ satisfy
\begin{equation}
|A_{kl}| <  c e^{-\gamma |k-l|^{\lambda}} \qquad\text{for}\,\,\,
\gam,\lam>0,
\notag 
\end{equation}
and some constant $c >0$. If $\lam=1$ we simply write $\expmata$.\\
(ii) $A$ belongs to the space $\polmat$ if the coefficients
$A_{kl}$ satisfy
\begin{equation}
|A_{kl}| <  c (1+|k-l|)^{-s} \qquad \text{for}\,\,\, s > 1,
\notag 
\end{equation}
and some constant $c>0$. 
\end{definition}

The following result is due to Jaffard~\cite{Jaf90}.

\begin{theorem}
\label{th:jaffard}
Let $A: \ltsp({\cal I}) \rightarrow \ltsp({\cal I})$ be an invertible
matrix, where ${\cal I}$ is ${\cal I} = \Zst,\Nst$ or
$\{0,\dots,N-1\}$. \\
(a) If $A \in \expmata$, then $\AI \in \expmatai$ for some $\gami <
\gam$.\\
(b) If $A \in \polmat$, then $\AI \in \polmat$.
\end{theorem}

For finite-dimensional matrices these results
(and in particular the involved constants) should be interpreted as follows. 
Think of the $n \times n$ matrix $\An$ as a finite section of an 
infinite-dimensional matrix $A$.
If we increase the dimension of $\An$ (and thus consequently the dimension
of $(A_n)^{-1}$) 
we can find uniform constants independent of $n$ such the
corresponding decay properties hold. This is of course not possible
for arbitrary finite-dimensional invertible matrices.

Theorem~\ref{th:jaffard}(a) shows that the entries of $\AI$ still
have exponential decay, however $\AI$ is in general not in the same
algebra as $A$, since we may have to use a smaller exponent. However in 
Theorem~\ref{th:jaffard}(b) both, the matrix $A$ and its inverse $\Ai$ 
belong to the same algebra, the quality of decay does not change.

From this point of view Theorem~\ref{th:jaffard}(b) is the most striking
result. The proof of Theorem~\ref{th:jaffard}(b) is rather delicate and 
lengthy. For biinfinite Toeplitz-type matrices 
this result can be proven much shorter (and extended to other types of 
decay) by using classical results of Gelfand et.~al.~on certain commutative 
Banach algebras. The following theorem is a weighted version of Wiener's
Lemma. It is implicitly contained in~\cite{GRS64}, but since it may be of 
independent interest we state and prove it explicitly.

\begin{theorem}
\label{th:gelfand}
Let $A=\{a_{kl}\}$ be a hermitian positive definite biinfinite Toeplitz matrix
with inverse $A^{-1}=\{\alpha_{kl}\}$. Let $v(k)$ be a positive (weight)
function with 
$$v(k+l)\le v(k)v(l)\,,$$
such that
\begin{equation}
\sum_{k=-\infty}^{\infty} |a_k| v(k) < \infty.
\label{gelfand1}
\end{equation}
If
\begin{equation}
\lim_{\ntoinf} \frac{1}{^n\!\sqrt{v(-n)}} = 1
\quad \text{and}\,\,\,
\lim_{\ntoinf}\, ^n\!\sqrt{v(n)} = 1 ,
\label{gelfandradius}
\end{equation}
then
\begin{equation}
\sum_{k=-\infty}^{\infty} |\alpha_k| v(k) < \infty.
\label{gelfand2}
\end{equation}
In particular, 
\begin{align}
\text{if}& \quad A \in \polmat \quad \text{for $s>1$, then} \,\,\,
\Ai \in \polmat ;
\label{gelfand3} \\
\text{if}& \quad A \in \expmatab \quad \text{for $0 < \lam < 1$, then}
\,\,\,
\Ai \in \expmatab .
\label{gelfand4}
\end{align}
\end{theorem}

\begin{proof}
Since $A$ is positive definite we have
\begin{equation}
f(\omega) = \sum_{k=-\infty}^{\infty} a_k e^{2\pi i k \omega} > 0
\label{e1}
\end{equation}
and by the properties of Laurent operators~\cite{GGK93}
$$1/f(\omega) = \sum_{k=-\infty}^{\infty} \alpha_k e^{2\pi i k \omega},
\qquad \text{where} \quad (A^{-1})_{k,l}=\alpha_{k-l}.$$ 
We denote by $W[v]$ the set of all formal series
$f=\sum_{k=-\infty}^{\infty}a_k X^k$ for which
\begin{equation}
\|f\|=\sum_{k=-\infty}^{\infty} |a_k| v(k) < \infty.
\notag 
\end{equation}
It follows from Chapter 19.4 of~\cite{GRS64} that $W[v]$ is a Banach algebra
with respect to the multiplication (discrete convolution)
\begin{equation}
fg = \sum_{l=-\infty}^{\infty}c_l X^l =
\sum_{l=-\infty}^{\infty}\left( \sum_{k=-\infty}^{\infty}a_{l-k}
b_k\right)X^l,
\notag 
\end{equation}
where $f=\sum_{k} a_k X^k$ and $g=\sum_{k}b_k X^k$.
By Theorem~2 on page~24 in~\cite{GRS64} an element of $W[v]$ has an
inverse in $W[v]$ if it is not contained in a maximal ideal of $W[v]$.
Any maximal ideal of $W[v]$ consists of elements of the form
(cf.~Chapter 19.4 in~\cite{GRS64})
\begin{equation}
\sum_{k=-\infty}^{\infty}a_k \xi^k = 0, \notag 
\end{equation}
where $\xi = \rho e^{2\pi i \omega}$ with
\begin{equation}
\rho_1 \le \rho \le \rho_2, \notag
\end{equation}
and 
\begin{equation}
\rho_1 = \underset{\ntoinf}\lim \frac{1}{^n\sqrt{v_{-n}}}
\,\,\,\,\text{and} \,\,\,
\rho_2 = \underset{\ntoinf}\lim ^n\sqrt{v_{n}}.
\notag
\end{equation}
Due to assumption~\eqref{gelfandradius} we get $\rho_1=\rho_2 =1$,
hence $\rho=1$. Thus a necessary and sufficient condition for an
element in $W[v]$ to be not contained in a maximal ideal of $W[v]$ 
is $\sum_{k}a_k e^{2\pi i k \omega} \neq 0$ for all $\omega$.
By assumption $A$ is positive definite, hence 
$f(\omega)=\sum_{k}a_k e^{2\pi i k \omega}>0$ for all $\omega$
and~\eqref{gelfand2} follows. 

Statements~\eqref{gelfand3} and~\eqref{gelfand4} are now clear, since in both 
cases we can easily find a weight function such that~\eqref{gelfand1} and
\eqref{gelfandradius} are satisfied.
\end{proof}

\begin{remark}
{\em
(i) Theorem~2.11 in~\cite{Dom56} by Domar and Theorem V B in~\cite{Beu38} 
by Beurling are closely related to Theorem~\ref{th:gelfand}. Their results are 
concerned with (non)quasi-analytic functions, for which they have to impose 
the more restrictive condition 
$$
\sum_{k=1}^{\infty}\frac{\log [v(kx)]}{k^2} < \infty,
\qquad \text{for all $x$,}
$$
on the weight function (called Beurling-Domar condition in~\cite{Rei68}).
For instance the function $v(k) = \exp(\frac{|k|}{1+\log(|k|)}), 
k \neq 0$ satisfies condition~\eqref{gelfandradius}, but not the Beurling-Domar 
condition. \\
(ii) Using Theorem~8.1 on page~830 in~\cite{GGK93} we can extend 
Theorem~\ref{th:gelfand}
to biinfinite block-Toeplitz matrices with finite-dimensional non-Toeplitz
blocks (i.e., Laurent operators with matrix-valued symbol). These
matrices play an important role in filter bank theory~\cite{Str98a}. \\
(iii) Note that $v(n)=\exp(\gamma n)$ does not satisfy 
condition~\eqref{gelfandradius},
that is why we have to introduce an exponent $\gamma_1 < \gamma$ in order
to estimate the decay of $\AI$, cf.~also Theorem~\ref{th:jaffard}.
However if $A \in \expmatab$ with $\lam <1$, then
condition~\ref{gelfandradius} is satisfied and -- as we have seen --
the decay of the entries of $\AI$ can be bounded by using the {\em same} 
parameters $\gamma,\lambda$.}
\label{remark1}
\end{remark}

\section{Approximation of infinite-dimensional Toeplitz-type systems} 
\label{s:infinite}

Infinite Toeplitz systems arise for instance in the discretization
of Wiener-Hopf integral equations or, more generally, in
one-sided infinite convolution equations, see~\cite{GF74}.
Biinfinite Toeplitz-type systems are encountered in doubly infinite
(discrete) convolution equations, as well as e.g.\ in filter bank 
theory~\cite{Str98a} or in the inverse heat problem~\cite{CN96}. 
In order to solve these problems we have to introduce a
finite-dimensional model.

For let $A:\ltZ \mapsto \ltZ$ be a hermitian positive definite ({\em hpd}
for short) biinfinite
Toeplitz matrix given by $\{a_{k,l}\}_{k=-\infty}^{\infty}$.
Let $y=\{y_k\}_{k=-\infty}^{\infty} \in \ltZ$. We want to solve
the system $Ax=y$. 

For $n \in \Nst$ and $y \in \ltZ$ define the orthogonal projections $\PN$ by
\begin{equation}
\label{defP}
\PN y = (\dots, 0,0, y_{-n+1},\dots , y_{n-1}, 0,0, \dots).
\end{equation}
By identifying the image of $\PN$ with the $2n-1$-dimensional space
$\Cst^{2n-1}$ we can express the $(2n-1)\times (2n-1)$ matrix $\An$ as
\begin{equation}
\notag 
\An = \PN A \PN,
\end{equation}
where we have used that $P^{\ast}=P$.
The $n$-th approximation $\xn$ to $x$ is then given by the
solution of the finite-dimensional system of equations
\begin{equation}
\notag 
\An \xn = \yn\,
\end{equation}
where $\yn :=\PN y$.

If $A$ is a singly infinite Toeplitz matrix and $y \in \ltsp(\Nst)$,
we proceed analogously by defining $\PN$ as
\begin{equation}
\notag 
\PN y = (y_{0}, y_{1},\dots , y_{n-1}, 0,0, \dots).
\end{equation}
This approach to approximate the solution of $Ax=y$ is usually called 
the {\em finite section method}, cf.~\cite{GF74}. 

The first question that arises when considering this method is 
``does $\xn$ converge to x?''. For the case when $A$ is not hpd this 
question has lead to deep mathematical results. See the book~\cite{GF74} 
and chapter~7 in~\cite{BS90} for more details.
For the case when $A$ is hpd the answer is easy and always 
positive. To see this, recall that since $A$ is hpd it follows that
$A_n$ is also hermitian positive definite, see~\cite{HJ94a}.
Furthermore, $\|A_n\| \le \|A\|$ and $\|(A_n)^{-1}\| \le \|A^{-1}\|$ 
for $n=1,2,\dots$.
Applying the Lemma of Kantorovich~\cite{RM94} yields that $(A_n)^{-1}$ 
converges strongly to $\AI$ for $\ntoinf$, i.e., $\xn$ converges to $x$ 
in the $\ltsp$-norm for any $y \in \ltZ$ (or for any $y\in \ltsp(\Nst)$
if $A$ is singly infinite).

An important aspect for applications is if we can give an estimate on how 
fast $\xn$ converges to $x$. It will be shown that the rate of 
approximation depends on the decay behavior and the condition number of the
matrix.

\begin{theorem}
\label{th1}
Let $Lx=y$ be given, where $L=\{a_{k,l}\}$ is a hermitian positive definite
biinfinite Toeplitz matrix and denote $\xn = \Lni \yn$.\\
(a) If there exist constants $c, c'$ such that
\begin{equation}
|a_k| \le c e^{-\gamma |k|} \enspace \text{and} \enspace
|y_k| \le c' e^{-\gamma |k|}, \,\,\,\gamma>0
\end{equation}
then there exists a $\gamma_1$ with $0 < \gamma_1 < \gamma$ and a constant 
$c_1$ depending only on $\gamma_1$ and on the condition number of $L$ such that
\begin{equation}
\|x-\xn\| \le c e^{-\gamma_1 n}.
\end{equation}
(b) If there exist constants $c, c'$ such that
\begin{equation}
|a_k| \le c(1+|k|)^{-s}  \enspace \text{and} \enspace
|y_k| \le c'(1+|k|)^{-s},\,\,\,s>1,
\end{equation}
then there exists a constant $c_1$ depending only on the condition 
number of $L$ such that
\begin{equation}
\label{polest1}
\|x-\xn\| \le c_1 n^{(1-2s)/2}.
\end{equation}
\end{theorem}
\begin{proof}
We have
\begin{align}
\|x-\xn\| & = \|\LI y - \Lni \yn \| \le
\|\LI\| \|y - L \Lni \yn \|  \notag \\
& \le \|\LI\| \big(\|y - \yn\| + \| (\Ln - L) \Lni \yn\|\big).
\label{mainin}
\end{align}

To prove statement~(a) we note that by Theorem~\ref{th:jaffard}(a) 
there exists a $\gamma_2 < \gamma$ such that 
$(\Lni)_{kl} \le c_2 e^{-\gamma_2 |k-l|}$ with a constant $c_2$ depending
only on $\gamma_2$ and on the condition number of $\Ln$. Since
$\sigma(\Ln) \subseteq [\fmin,\fmax]$ we get $\cond(\Ln) \le \cond(L)$ for 
all $n$. That means we can choose $c_2$ independently of $n$.
Write $\zn = (\Ln - L) \Lni \yn$ and note that $\zn_k = 0$
for $|k| < n$. Since the non-zero entries of $(\Ln-L)$ decay
exponentially, it is easy to show that there exists a $\gamma_1$ with
$0 <\gamma_1 <\gamma_2$ such that $\|\zn\| \le c_3 e^{-\gamma_1 n}$
for some constant $c_3$.  It is trivial that $\|y-\yn\|$ also decays 
exponentially for $\ntoinf$. We absorb $\|\LI\|$ and the other constants
in the constant $c_1$ and get the desired result.

For the proof of part (b) we proceed analogously to above by applying 
Theorem~\ref{th:jaffard}(b) to conclude that 
$$|(\Lni \yn)_k| \le c_2 (1+k)^{-s}$$
for some constant $c_2$ depending only on $\cond(L)$ and on $s$.
The norm $\|y-\yn\|$ can be estimated via
\begin{gather}
\|y-\yn\|^2=\sum_{|k|\ge n} |y_k|^2 \le 
2c \sum_{k=n}^{\infty}(1+k)^{-2s} \le 2c \int \limits_{n-1}^{\infty} 
(1+x)^{-2s} dx \le 2c \frac{n^{1-2s}}{2s-1},
\label{e13}
\end{gather}
similarly for $\|\zn\|$ where $\zn:=(\Ln-L) \Lni\yn$.
Since all arising constants - absorbed in one constant $c_1$ -
depend only on $\cond(L)$ and on the exponent $s$, the proof is complete.
\end{proof}
\begin{remark} 
\label{remark2}
{\em
Theorem~\ref{th1} holds if we replace the system $Lx=y$ by 
a singly infinite Toeplitz system $Tx=y$ with corresponding decay
conditions on $T$ and $y$ and approximate its solution by considering
the finite system $\Tn \xn =\yn$.
}
\end{remark} 

The proof of Theorem~\ref{th1} is essentially based on the fact that under 
appropriate decay conditions on $L$ and $y$ (resp.\  $\Ln$ and $\yn$) 
$\Li$ and $\Lni \yn$ have similar decay properties. Thus, if one can
show that $\Lni$ has the same decay properties as $\Li$ one can use
Theorem~\ref{th:gelfand} in order to generalize Theorem~\ref{th1} to various
other decay conditions. This may however not always
lead to simple and closed-form expressions for $\|(\Ln-L) \Lni \yn\|$, 
therefore I have restricted myself to the most frequently encountered 
decay properties.

{\sc Example 1:}
We illustrate Theorem~\ref{th1} by a numerical example. We consider
$Lx=y$, where $L$ is the biinfinite Toeplitz matrix with entries
$a_k = (1+|k|)^{-s}, k \in \Zst$ for $s=2$ and $y$ consists of
random entries having the same polynomial decay rate as the entries $a_k$. 
To compare the error $\|x-\xn\|$ with the error estimate~\eqref{polest1} 
we would need the true solution $x$. Since the solution of this 
biinfinite system cannot be computed analytically we compute the 
``true'' solution of $Lx=y$ by solving $L_{n_0} x^{(n_0)} = y^{(n_0)}$ for 
very large $n_0$ (we choose $n_0=32768$). Using~\eqref{mainin} 
and~\eqref{condest6} we can estimate that in the worst case 
$\|x-x^{(32768)}\| \approx 10^{-6}$, so that $x^{(32768)}$ can mimick the 
true solution with sufficiently high  accuracy for this experiment.

\begin{figure}
\begin{center}
\epsfig{file=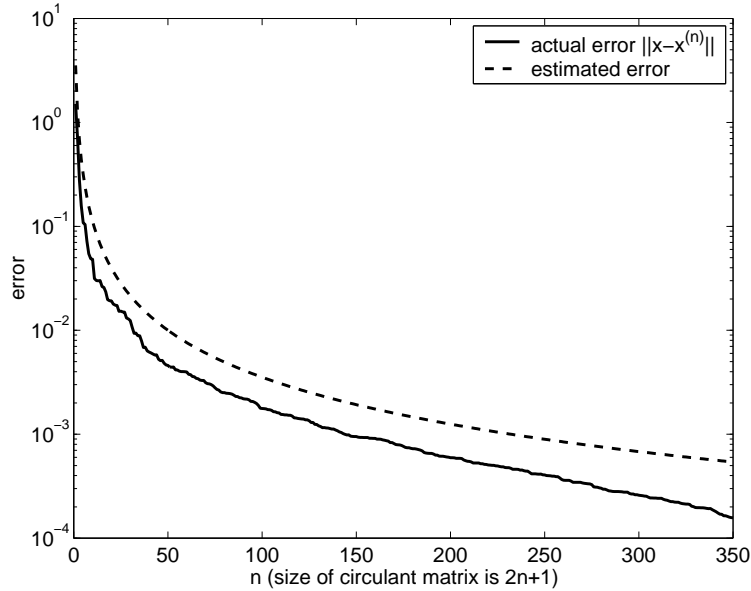,width=100mm}
\caption{Actual approximation error and error estimate~\eqref{polest1} 
from Theorem~\ref{th1}(b) for the system $Lx=y$, where $L$ is a biinfinite 
hermitian Toeplitz matrix with polynomial decay.}
\label{fig:fig1}
\end{center}
\end{figure}

Then we approximate this solution by the finite section method
as in Theorem~\ref{th1} for $n=0,\dots, 350$ and compute for each $n$
the error $\|x-\xn\|$  as well as the error 
estimate in~\eqref{polest1}. Note that an explicit expression for
the constant $c_1$ in~\eqref{polest1} is not known, we only know 
that it depends on the condition number of $L$.
In this example we use $c_1=\cond(L)$ (a different example may 
require a different choice). The result, illustrated in
Figure~\ref{fig:fig1}, shows that the asymptotic behavior of
the error rate is well estimated by the given error bound.

It is well-known that the product of two Laurent operators and the
inverse of a Laurent operator (if it exists) is again a Laurent
operator. This is of course not true for singly infinite or finite
Toeplitz matrices (and this is one of the reasons which makes 
the ``Toeplitz business'' so interesting). Hence one may argue that the
``canonical'' finite-dimensional analogue of Laurent operators are
not Toeplitz matrices but circulant matrices, since they also form
an algebra. 
Thus for a given biinfinite hermitian Toeplitz matrix $L$ with entries
$L_{kl}=a_{k-l}$ we define the hermitian circulant matrix $\Cn$ of size
$(2n-1)\times (2n-1)$ by
\begin{equation}
\Cn = 
\begin{bmatrix}
a_0      & \conj{a}_1 & \dots       & \conj{a}_{n-2} & \conj{a}_{n-1} & 
a_{n-1}  &  a_{n-2}   & \dots       & a_1 \\
a_1      &  a_0       & \conj{a}_1  & \dots          &         &
               &                &       & a_2 \\
\vdots         &                &\ddots &         &         &
               &                &       &  \vdots    \\   
\conj{a}_1            &                &\dots  &         &         &
               &                &a_1 & a_0    
\end{bmatrix}.
\label{circmat}
\end{equation}
We also say that $\Cn$ is generated by $\{a_k\}_{k=-n+1}^{n-1}$.

\begin{remark}
{\em
$\Cn$ does not have to be positive definite if $L$ is positive definite,
e.g.~see~\cite{CS89}. However - as pointed out in~\cite{CS89} - if
$L$ is at least in Wiener's algebra then one can always find an $N$ such
that $\Cn$ is invertible for all $n>N$. The faster the decay of the entries
of $L$ the smaller this $N$ has to be.}
\label{remark3}
\end{remark}

We can do even a little better and
estimate how well the extrema of the defining function of $L$ are
approximated by the extreme eigenvalues of $\Cn$.

\begin{lemma}
\label{le:cond}
Let $L$ be a biinfinite hermitian Toeplitz matrix with entries
$L_{k,l}=a_{k-l}$ where $a=\{a_k\}_{k=-\infty}^{\infty}$ and set 
$f(\omega)=\sum_{k=-\infty}^{\infty} a_k e^{2 \pi i k \omega}$. 
Let $C_m$ be the associated $(2m-1) \times (2m-1)$
circulant matrix with first row 
$(a_0, \conj{a}_1, \dots, \conj{a}_{m-1}, a_{m-1}, \dots, a_1)$.
Denote the maximum and minimum eigenvalue resp.\ of $C_m$ by 
$\lambda^{(m)}_{\max}$ and $\lambda^{(m)}_{\min}$. \\
(a) If $L$ is $n$-banded with $m>n$, then
\begin{gather}
\lambda^{(m)}_{\max} \le \fmax \le \lambda^{(m)}_{\max} + 
2 \sin \Big(\frac{\pi n}{2(2m-1)^2}\Big)\|a\|_1, \label{condest1}\\
\lambda^{(m)}_{\min} \ge \fmin \ge \lambda^{(m)}_{\min} - 
2 \sin \Big(\frac{\pi n}{2(2m-1)^2}\Big)\|a\|_1. \label{condest2}
\end{gather}
(b) If $|a_{k}| \le c e^{-\gamma |k|}$ for $k \in \Zst, c>0$, then
\begin{gather}
|\fmax - \lambda^{(m)}_{\max}| \le \frac{2c}{1-e^{-\gamma}}
\Big[2\sin \Big(\frac{\pi m}{2(2m-1)^2}\Big) + e^{-\gamma m}\Big],
\label{condest4}
\end{gather}
a similar estimate holds for $|\fmin - \lambda^{(m)}_{\min}|$.\\
(c) If $|a_{k}| \le c (1+|k|)^{-s}$ for $k \in \Zst, s>1, c>0$, then
\begin{gather}
|\fmax - \lambda^{(m)}_{\max}| \le \frac{2c}{s-1} 
\Big[2\sin \Big(\frac{\pi m}{2(2m-1)^2}\Big) + m^{1-s}\Big],
\label{condest6}
\end{gather}
a similar estimate holds for $|\fmin - \lambda^{(m)}_{\min}|$.
\end{lemma}
\begin{proof}
(a): It is well-known~\cite{Dav79} that the eigenvalues of $C_m$ are given by 
\begin{equation}
\sum_{k=-m+1}^{m-1} a_k e^{2\pi i k l/(2m-1)}, \qquad 
\text{for $l=-m+1,\dots,m-1$.}
\label{e18}
\end{equation}
For case (a) that means 
they are regularly spaced samples $f(\frac{l}{2m-1})$ of the function $f$, 
which immediately yields the left hand side of the 
inequalities~\eqref{condest1} and~\eqref{condest2}.
In order to prove the right hand side of~\eqref{condest1}
and~\eqref{condest2} it is sufficient to estimate
\begin{equation}
\underset{\omega,l}{\max} \Bigl|f(\frac{l+\omega}{2m-1})-f(\frac{l}{2m-1})\Bigr|,
\notag 
\end{equation}
where $\omega \in [-\frac{1}{2(2m-1)},\frac{1}{2(2m-1)}]$ and $l=-m+1,\dots,m-1$.
Define the sequence $\{\tilde{a}_k\}_{k=-m}^m$
by $\tilde{a}_k=a_k$ if $|k| \le n$ and $\tilde{a}_k=0$ if $|k| > n$.
There holds
\begin{gather}
\underset{\omega,l}{\max} \Bigl|f(\frac{l}{2m-1})-f(\frac{l+\omega}{2m-1})\Bigr| 
= \underset{\omega,l}{\max} \Bigl|\sum_{k=-m}^{m} \tilde{a}_k e^{2\pi i k l/(2m-1)}
(e^{2\pi i k \omega/(2m-1)} -1)\Bigr| \notag \\
\le \underset{\omega,|k| \le n}{\max}\bigl|e^{2\pi i k\omega/(2m-1)}-1\bigr|
\sum_{k=-n}^{n} |a_k|  \notag \\
\le
\underset{\omega,|k| \le n} {\max} 2 |\sin (\pi k \omega/(2m-1))| \|a\|_1
\le 2 \sin \big(\frac{\pi m}{2(2m-1)^2}\big) \|a\|_1.
\label{e4}
\end{gather}
Relations~\eqref{condest1} and \eqref{condest2} follow now from this
estimate.

Statements (b) and (c) can be proved similarly by using 
\begin{gather}
|f\big(\frac{l+\omega}{2m-1}\big) - \lambda_l^{(m)}| \le
\sum_{k=-m+1}^{m-1}|a_k| |e^{2\pi i k\omega/(2m-1)}-1| + \sum_{|k|\ge m} |a_k|,
\label{e14}
\end{gather}
and applying the corresponding decay properties to~\eqref{e14}.
\end{proof}

\begin{remark}
\label{remark4}
{\em
The left part of inequality~\eqref{condest2} reads $\lambda^m_{\min} \ge \fmin$.
This implies that Strang's preconditioner is always positive definite
for $n$-banded hermitian Toeplitz matrices of size $> (2n \times 2n)$
with positive generating function. Hence the ``sufficiently large
$n$''-condition at the end of section~2 in~\cite{HN94} can be omitted.}
\end{remark}

It is obvious that decay properties for circulant matrices cannot
be defined in the same way as for non-circulant matrices. Hence, by stating 
that $C_m$ has, say, exponentially decaying entries, we mean that
the generating sequence $\{a_k\}_{k=-m+1}^{m-1}$ satisfies
$|a_k|\le c e^{-\gamma |k|}$, in which case the entries of $C_m$ will
decay exponentially off the corners of the matrix (instead of off
the diagonal).
In analogy to the theorems in section~\ref{s:decay} it is natural to ask if the 
inverse of $C_m$ also inherits these decay properties.
The following theorem shows that at least for $m$ sufficiently large
this is the case. The entries of $C_m^{-1}$ uniformly approximate the entries
of $\Li$ with an error rate depends on the decay properties and the
condition number of $L$.

\begin{theorem}
\label{th:circulant}
Let $L$ be a hermitian positive definite Laurent operator with entries
$L_{kl}=a_{k-l}$ and let $C_m$ be the
associated circulant $(2m-1)\times (2m-1)$ matrix as defined in
\eqref{circmat}. Denote the entries of $L^{-1}$ by 
$(L^{-1})_{kl}=\{\alpha_{k-l}\}$ and let $(C_m)^{-1}$ (if it exists)
be generated by $\{\beta_k\}_{k=-m+1}^{m-1}$.\\
(i) If $L \in \expmata$, then for sufficiently large $m$
\begin{equation}
|\alpha_k - \beta_k| \le  c e^{-\gamma_1 m},
\label{mainexp}
\end{equation}
with $0 <\gamma_1 \le \gamma$ and some constant $c$ depending on
$\cond(L)$ and $\gamma_1$. \\
(ii) If $A \in \polmat$ for $s>1$, then for sufficiently large $m$
\begin{equation}
\label{mainpol}
|\alpha_k - \beta_k| \le  c \frac{m^{1-s}}{s-1}
\end{equation}
with some constant $c$ depending on $\cond (L)$.
\end{theorem}

\begin{proof}
(i): Set $f(x)=\sum_{k=-\infty}^{\infty} a_k e^{2\pi i k x}$ and
$f_m(x)=\sum_{k=-m+1}^{m-1} a_k e^{2\pi i k x}$. By the properties of
Laurent operators~\cite{GGK93} $\{\alpha_k\}_{k \in \Zst}$ is given by
\begin{equation}
\alpha_k = \int \limits_{-1/2}^{1/2}\finvx \expkx dx , 
\qquad k\in \Zst,
\label{e5}
\end{equation}
By Remark~\ref{remark3} and Lemma~\ref{le:cond} we can easily find an $N$ 
such that $C_m$ is invertible for all $m>N$, which implies
that $f_m > 0$. In this case by the properties of circulant 
matrices~\cite{Dav79} the entries $\{\beta_k\}_{k=-m+1}^{m-1}$ of
$C_m^{-1}$ can be computed as
\begin{equation}
\beta_k = \frac{1}{2m-1}\sum_{l=-m+1}^{m-1}\fminvlm \explm  , 
\qquad k=-m+1,\dots,m-1.  \label{e6}
\end{equation}
Now consider
\begin{gather}
|\alpha_k - \beta_k| = 
\Bigl|\int \limits_{-1/2}^{1/2} \finvx \expkx dx-
\sum_{l=-m+1}^{m-1} \fminvlm \explm \Bigr| \notag \\
\le \Bigl|\int \limits_{-1/2}^{1/2} \finvx \expkx dx - 
\frac{1}{2m-1}\sum_{l=-m+1}^{m-1} \finvlm \explm \Bigr| + \label{eq3} \\
+ \frac{1}{2m-1}\Bigl|\sum_{l=-m+1}^{m-1} \finvlm \explm - 
\sum_{l=-m+1}^{m-1} \fminvlm \explm \Bigr| 
\label{eq4}
\end{gather}
We estimate the expression above in two steps:\\
1. We first consider~\eqref{eq3}. Note that 
$\finvx = \sum_{k=-\infty}^{\infty} \alpha_k \expkx$ and
\begin{equation}
\finvlm = \sum_{k=-\infty}^{\infty} \alpha_k \explm
= \sum_{p=-m+1}^{m-1} \left(\sum_{q=-\infty}^{\infty} \alpha_{p+(2m-1)q}\right)
e^{2\pi ilp/(2m-1)}.
\notag
\end{equation}
By setting
$\delta_p = \sum_{q=-\infty}^{\infty}\alpha_{p+(2m-1)q}$
we get
$$\finvlm = \sum_{p=-m+1}^{m-1}\delta_p e^{2\pi i l p/(2m-1)},
\qquad l=-m+1,\dots,m-1,$$
and
$$\delta_k=\frac{1}{2m-1}\sum_{p=-m+1}^{m-1}\finvpm e^{-2\pi i kp/(2m-1)}.$$
Hence
\begin{gather}
\Bigl|\int \limits_{-1/2}^{1/2} \finvx \expkx dx - 
\frac{1}{2m-1}\sum_{l=-m+1}^{m-1} \finvlm \explm \Bigr | \notag \\
= | \alpha_k - \delta_k | \le \sum_{q \neq 0} |\alpha_{k+(2m-1)q}|.
\label{e11}
\end{gather}
By Theorem~\ref{th:jaffard}(a) $A \in \expmata$ implies $\AI \in \expmatai$.
Hence we get for $k=-m+1, \dots, m-1$
\begin{gather}
\sum_{q\neq 0}|\alpha_{k+(2m-1)q}| \le 
2c_1 \sum_{q=1}^{\infty} e^{-\gamma_1 (k+(2m-1)q)} \notag \\
\le 2c_1 \sum_{q=1}^{\infty} e^{-\gamma_1 ((2m-1)q-m+1)}
= \frac{2c_1 e^{-\gamma_1 m}}{1-e^{-\gamma_1(2m-1)}}. \label{e15}
\end{gather}
2. Now we estimate~\eqref{eq4}:
\begin{gather}
\frac{1}{2m-1}\Bigl|\sum_{l=-m+1}^{m-1} \finvlm \explm -
\sum_{l=-m+1}^{m-1} \fminvlm \explm \Bigr| \notag \\ 
\le\frac{1}{2m-1} \sum_{l=-m+1}^{m-1}\bigl|\finvlm - \fminvlm \bigr| \notag \\
\le\frac{1}{2m-1} \sum_{l=-m+1}^{m-1} \bigl|\finvlm  \bigr|  \bigl|\fminvlm \bigr|
\bigl|\flm -\fmlm \bigr| \notag \\
\le\frac{1}{2m-1} \underset{|l|\le m}{\max} \bigl|\finvlm \bigr| \bigl|\fminvlm \bigr|
\sum_{l=-m+1}^{m-1}\sum_{|k|>m}|a_k| \notag \\
\le \|A^{-1}\|\|C_m^{-1}\|\sum_{|k|>m}|a_k|.\label{e12}
\end{gather}
By Lemma~\ref{le:cond} we can easily find for any $\eps >0$ an $N$ such that 
for all $m >N$ there holds $\|C_m^{-1}\| \le (1+ \eps)\|A^{-1}\|$. Thus
\begin{gather}
\|A^{-1}\|\|C_m^{-1}\| \sum_{|k|>m}|a_k| 
\le (1+\eps)\|A^{-1}\|^2 2 \sum_{k=m}^{\infty}|a_k| \notag \\
\le (1+\eps)\|A^{-1}\|^2 2 \frac{e^{-\gamma m}}{1-e^{-\gamma}}
\le (1+\eps)\|A^{-1}\|^2 c_2 e^{-\gamma_1 m}
\label{eq5}
\end{gather}
for some $\gamma_1 < \gamma$ and some constant $c_2$.
By combining~\eqref{e11}, \eqref{e15}, and \eqref{eq5} and 
hiding expressions as $(1+ \eps)\|A^{-1}\|$ in the constant $c$, we
obtain estimate~\eqref{mainexp}.

(ii): The proof of~\eqref{mainpol} is similar to the
proof of~\eqref{mainexp}. The only steps that require a modification are
\eqref{e15} and~\eqref{eq5}.
By Theorem~\ref{th:jaffard}(a) $A \in \polmat$ implies $\AI \in \polmat$.
Hence we can estimate $\sum_{q \neq 0} |\alpha_{k+(2m-1)q}|$ 
as follows.
\begin{gather}
\sum_{q \neq 0} |\alpha_{k+(2m-1)q}| \le
c' \sum_{q \neq 0} (1+|(2m-1)q+k|)^{-s} \notag \\
= 2c'(1+|2m-1+k|)^{-s} + 2c' \sum_{q=2}^{\infty}(1+|(2m-1)q+k|)^{-s}\notag \\
\le 2c'(1+|2m-1-m+1|)^{-s} + 2c'
\sum_{q=2}^{\infty}(1+|(2m-1)q-m+1|)^{-s}\notag \\
\le 2c' (1+m)^{-s} + 2c'\int \limits_{1}^{\infty}(1+|(2m-1)x-m+1|)^{-s}\,dx
\notag \\
=2c'(1+m)^{-s} + 2c' \frac{(1+m)^{-s}}{(2m-1)(s-1)}
\le \frac{2c (1+m)^{-s}}{s-1}.
\label{e16}
\end{gather}
By adapting~\eqref{eq5} to the case of polynomial decay we can estimate
\eqref{eq4} by 
\begin{equation}
\label{e17}
2c\|A^{-1}\| \|C_m^{-1}\| \frac{m^{1-s}}{s-1}.
\end{equation}
Combining~\eqref{e16} with~\eqref{e17} yields the desired result.
\end{proof}

The main portion of the proof of statement (ii) above is
due to Gabriele Steidl~\cite{Ste00}.

\section{Error estimates for approximate solution of deconvolution problems}
\label{s:deconv}

Consider the convolution of two sequences $a=\{a_k\}_{k=-\infty}^{\infty}, 
c =\{c_k\}_{k=-\infty}^{\infty} \in \ltZ$, given by $b = a \ast c$. Here
$a$ may represent an impulse response or a blurring function.
Given $a$ and $b$ our goal is to compute $c$. This is known as deconvolution. 
In matrix notation the problem can be expressed as $Lc=b$ where $L$ is a
biinfinite Toeplitz matrix with entries $L_{kl}=a_{k-l}$. 

Sometimes $a$ and $c$ have compact support, in which
case $b$ also has compact (although larger) support and $c$ can be
computed by solving a finite banded Toeplitz system, see~\cite{NPT96}.
It is well-known that this can be done efficiently via FFT by embedding 
the Toeplitz matrix into a circulant matrix. Of course, this approach is 
very attractive from a numerical viewpoint, at least if the system
is well-conditioned (we discuss the ill-conditioned case in 
section~\ref{s:precond}).

However, if either $a$ or $c$ does not have compact support, the
reduction of $Lc=b$ to a finite Toeplitz system obviously
will introduce a truncation error and embedding the Toeplitz matrix
into a circulant yields an additional (perturbation) error.
We nevertheless can try to make use of the FFT-based approach with the
hope to get a good approximation to the solution.
According to~\cite{GHK94} it has been shown in~\cite{Lev80} 
that the solution of doubly infinite convolution systems can be 
approximated by solutions of finite circulant systems. 
(Note that the approximation by finite circulant systems does not
apply to one-sided infinite convolution equations.) 

Hence we are concerned with the problem of how good the approximation is
obtained in that way and how fast the approximation converges to the 
true solution.

To answer these questions we proceed as follows. In the first step we
study the approximate solution of a finite Toeplitz system by circulant
embedding. In the second step we combine the obtained results with 
Theorem~\ref{th1} to derive estimates for the approximate solution of doubly 
infinite convolution equations by finite circulant systems. The results that 
we will collect in the first step will also be very useful in
section~\ref{s:precond} for the analysis of circulant preconditioners.

Let $\An \xn =\yn$ be given, where $\An$ is an hpd $n \times n$ Toeplitz
matrix. Of course, we have in mind that $\An$ and $\yn$ are finite sections 
of the biinfinite Toeplitz matrix $L$ and the right-hand side $y$,
respectively.

As usual, we embed $\An$ into a circulant matrix $\Cnn$ of size $2n \times 2n$
as follows:
\begin{equation}
\Cnn = 
\begin{bmatrix} 
\An & \Bna \\
\Bn & \An
\end{bmatrix},
\label{embed}
\end{equation}
where $\Bn$ is the $n \times n$ Toeplitz matrix with 
first row given by $b_k = a_{n-k}$ for $k=1,\dots,n-1$. If 
$a_n$ is known we set $b_0 = a_n$, otherwise we define $b_0=0$.
For the following considerations it does not matter if we embed $\An$
into a circulant matrix of size $(2n-1) \times (2n-1)$ or of size
$2n \times 2n$ (or larger). Choosing $\Bn$ to be of the same size
as $\An$ is just more convenient for the proofs below.
Since $\Cnn$ is circulant we can find $\An$ again in the center of
$\Cnn$ (and at any other position along the main diagonal of $\Cnn$),
and express the embedding of $\An$ into $\Cnn$ as follows
\begin{equation}
\Cnn = 
\begin{bmatrix} 
\times & \times & \times \\
\times &   \An  & \times \\
\times & \times & \times \\
\end{bmatrix}. 
\notag 
\end{equation}

In spite of the biinifite system $Lx=y$ in the background it is 
useful to embed $\yn$ symmetrically into a vector $\ynn$ of length $2n$ as 
follows 
\begin{equation}
\ynn = [0,\dots , 0, \yn, 0, \dots, 0].
\notag 
\end{equation}
We assume for the moment that $\Cnn$ is invertible (and will justify this
assumption later).
An approximate solution $\zn$ to $\An \xn = \yn$ is now obtained by
solving $\Cnn \znn = \ynn$ and setting $\zn =
\{\znn_k\}_{k=n/2+1}^{2n-n/2}$ (i.e., we take as approximation the central 
part of $\znn$ corresponding to the embedding of $\yn$ into $\ynn$).

We partition $\Cnni$ in the same way as $\Cnn$ as follows
\begin{equation}
\Cnni = 
\begin{bmatrix} 
M_n & T_n^{\ast} \\
T_n & M_n
\end{bmatrix}, 
\notag 
\end{equation}
where $M_n$ is a Toeplitz matrix and by Cauchy's interlace theorem
(cf.~\cite{HJ94a}) invertible. Define $\Sn:=M_n^{-1}$, then it is easy to 
see that $\zn$ can be obtained as the solution of 
$$\Sn \zn = \yn \,.$$

The question is now how well does $\zn$ approximate $\xn$.

\begin{theorem}
\label{th:circ2}
Let $\An \xn =\yn$ be given where $\An$ is an $n \times n$ hermitian positive 
definite Toeplitz matrix with $(\An)_{kl}=a_{k-l}$ and 
$a_k=\int \limits_{-\frac{1}{2}}^{\frac{1}{2}} 
f(\omega) e^{2\pi i \omega k} \,d\omega$.
Assume that $|a_k| \le c e^{-\gam |k|}$ and
$|(\yn)_k| \le c' e^{-\gam |n/2 - k|}$.
Suppose that $\Cnn$ as defined in~\eqref{embed} is invertible and
let $\zn$ be the solution of $\Sn \zn = \yn$, where $\Sn$ is the 
$n \times n$ leading principal submatrix of $\Cnni$.
Then there exists a $\gamma_1$ with $0 < \gamma_1 < \gamma$ and a constant 
$c_1$ depending only on $\fmin$ and $\fmax$ and on $\gamma_1$ such that
\begin{equation}
\|\xn-\zn\| \le c_1 e^{-\gamma_1 n} .
\label{e7}
\end{equation}
\end{theorem}

\begin{proof}
Similar to Theorem~\ref{th1} we write
\begin{equation}
\|\xn-\zn\| \le \|\Ani\| \|(\An - \Sn) \Sni \yn\| . 
\label{basic}
\end{equation}
Using the Schur complement~\cite{HJ94a} we can write $\Sn$ as
\begin{equation}
\Sn = \An - \Bn \Ani \Bna.
\label{schur}
\end{equation}
We will first show that the entries of the matrix $\An-\Sn$ are
exponentially decaying off the corners of the matrix. Note that
equation~\eqref{schur} implies 
\begin{equation}
\An - \Sn = \En
\notag 
\end{equation}
where $\En:= \Bn \Ani \Bna$.
We analyze the decay behavior of $\En$ in two steps by considering
first $\Bn \Ani$ and then $(\Bn \Ani) \Bna$.

Recall that $|\Ankl| \le c \egamkl$ and note that
$|\Bnakl| \le c \egamnkl$. By Theorem~\ref{th:jaffard} we know that
$|\Anikl| \le c_2 \egamiikl$ for some $0 <\gamii < \gam$, where $c_2$ depends 
on $\gamii$ and on $\fmin$ and $\fmax$, but is independent of $n$. 
We set $\delta = \gam - \gamii$.  There holds 
\begin{equation}
|(\Bn \Ani)_{kl}| = |\sum_{j=0}^{n-1} (\Bn)_{kj} (\Ani)_{jl} |
\le  c c_2 \sum_{j=0}^{n-1} \egamnkj \egamiijl .
\label{sum1}
\end{equation}
For simplicity we will absorb any constants arising 
throughout this proof that depend solely on $\gam$ (or $\gamii$) in the 
constant $c$.
We analyze the sum \eqref{sum1} further by splitting it up into three parts and
in addition consider first the entries $(\Bn \Ani)_{kl}$ with $k\ge l$.\\
(i) $0 \le j < l$: 
\begin{gather}
c \sum_{j=0}^{l-1} \egamnkj \egamiijl = 
c e^{-\gam (n-k)} e^{-\gamii l} \sum_{j=0}^{l-1} e^{-(\gam-\gamii)j} 
\le c e^{-\gamii(n-k+l)}.
\notag 
\end{gather}
(ii) $ l \le j < k$:
\begin{gather}
c \sum_{j=l}^{k-1} \egamnkj \egamiijl \le
c e^{-\gam (n-k)} e^{\gamii l}
\frac{e^{-(\gam+\gamii)l}}{1-e^{-(\gam+\gamii)}}
\le c e^{-\gamii(n-k+l)}.
\notag 
\end{gather}
(iii) $ k \le j < n$:
\begin{gather}
c \sum_{j=k}^{n-1} \egamnkj \egamiijl \le 
c e^{-\gam (n+k)} e^{\gamii l}
\frac{e^{(\gam-\gamii)(n-1)}}{1-e^{-(\gam-\gamii)}}\\
= c e^{-\gamii (n+k-l)} 
e^{-\delta (n+k)}\frac{e^{\delta(n-1)}}{1-e^{-(\gam-\gamii)}}
\le c e^{-\gamii(n-l+k)} \le c e^{-\gamii(n-k+l)}.
\notag 
\end{gather}
Similar expressions can be obtained for the case $l \ge k$ by interchanging
the roles of $k$ and $l$ in the derivations above. Thus
\begin{equation}
|(\Bn \Ani)_{kl}|\le c e^{-\gamii(n-|k-l|)}.
\notag 
\end{equation}

We now estimate the decay of the entries of $\Bn \Ani \Bn$.
Since $\Bn \Ani \Bna$ is hermitian, it is sufficient to consider
only the entries $(\Bn \Ani \Bna)_{kl}$ with $k\ge l$. 
There holds
\begin{equation}
|(\Bn \Ani \Bna)_{kl}| \le 
c \sum_{j=0}^{n-1} e^{-\gamii(n-|k-j|)} e^{-\gam (n-|j-l|)}.
\label{e8}
\end{equation}
As before we proceed by splitting up this sum into three parts.\\
(i) $0 \le j < l$:
\begin{gather}
c \sum_{j=0}^{l-1} e^{-\gamii(n-|k-j|)} e^{-\gam (n-|j-l|)}
=c e^{-(\gam+\gamii)n} e^{\gamii k} e^{\gam l} \frac{1}{1-e^{-(\gam+\gamii)}}
\le c e^{-\gamii(2n-k-l)} 
\notag
\label{case2a}
\end{gather}
(ii) $l \le j < k$:
\begin{gather}
c \sum_{j=l}^{k-1} e^{-\gamii(n-|k-j|)} e^{-\gam (n-|j-l|)}\le
c e^{-(\gamii+\delta)(n+l)}e^{-\gamii (n-k)}\frac{e^{\delta (k-1)}}{1-e^{-\delta}}
\le c e^{-\gamii(2n-k+l)} 
\notag
\label{case2b}
\end{gather}
(iii) $k \le j < n$:
\begin{gather}
c \sum_{j=k}^{n-1} e^{-\gamii(n-|k-j|)} e^{-\gam (n-|j-l|)}\le
e^{-(\gam+\gamii)n} e^{-\gamii(k+l)} 
\frac{e^{(\gam+\gamii)(n-1)}-e^{(\gam+\gamii)(k-1)}}{1-e^{-(\gam+\gamii)}}
\notag \\
\le c e^{-\gamii(k+l)}.
\notag
\label{case2c}
\end{gather}
Hence, by combining~(i), (ii), and (iii) we get
\begin{equation}
\label{bab}
|(\Bn \Ani \Bna)_{kl}| \le 
c (e^{-2\gamii n} e^{\gamii (k+l)} +  e^{-2\gamii n} e^{\gamii |k-l|)}
+  e^{-\gamii(k+l)}).
\end{equation}

The entries of $\Cnn$ satisfy
\begin{equation}
|(\Cnn)_{lk}|\le 
\begin{cases} 
c e^{-\gam |k-l|} & \text{for}\,\,|k-l|=0,\dots, n-1 ,\\
c e^{-\gam (2n-|k-l|)} & \text{for}\,\,|k-l|=n,\dots, 2n-1 .\\
\end{cases} 
\end{equation}
By Theorem~\ref{th:circulant} there exists a
$\gamma_3 < \gam$ and a constant $c_3$ depending on $\gamma_3$ and
on $\fmin$ and $\fmax$ such that
\begin{equation}
|(\Cnni)_{lk}|\le 
\begin{cases} 
c_3 e^{-\gamma_3 |k-l|} & \text{for}\,\,|k-l|=0,\dots, n-1 ,\\
c_3 e^{-\gamma_3 (2n-|k-l|)} & \text{for}\,\,|k-l|=n,\dots, 2n-1 .\\
\end{cases} 
\end{equation}
Hence $\zn:= \Sni \yn$ satisfies $|\zn_k| \le c_3 e^{-\gamma_3 |n/2-k|}$.
Set $\errn = (\An-\Sn) \zn$. After some lengthy but straightforward 
computations we get
\begin{equation}
|(\errn)_k| \le
\begin{cases}
c_1 e^{-\gamma_1 (n/2+k)} &\text{for $k=0,\dots,n/2$},\\
c_1 e^{-\gamma_1 (3n/2-k)} &\text{for $k=n/2,\dots,n$}
\end{cases}
\end{equation}
for some $\gamma_1 < \gamma_3$. Hence $\|\errn\| \le c e^{-\gamma_3 n/2}$,
which together with~\eqref{basic} completes the proof.
\end{proof}

\begin{corollary}
\label{corcirc}
Let $Lx=y$ be given where $L$ is a biinfinite hermitian positive definite
Toeplitz matrix with entries $L_{kl}=a_{k-l}$ and let $S_n$ be as defined 
in~\eqref{schur}. Assume $L \in \expmata$ and $|y_k| \le e^{-\gamma |k|}$
and let $\zn$ be the solution of $\Sn \zn = \yn$.
Then there exists an $N$ such that for all $n > N$
\begin{equation}
\notag 
\|x-\zn\| \le c e^{-\gamma_1 n},
\end{equation}
for some $0 < \gamma_1 < \gamma$ and a constant $c$ independent of $n$.
\end{corollary}
\begin{proof}
First note that by Remark~\ref{remark3} we can always find an $N$ such
that $\Sn$ exists for all $n >N$.
There holds
\begin{equation}
\|x-\zn\| = \|\Li y - \Sni \yn\| \le
\|\Li y - \Lni \yn\| + \|\Lni \yn - \Sni \yn\|
\label{e9}
\end{equation}
where $\An$ is an $n \times n$ finite section of $L$.
The result follows now by applying Theorem~\ref{th1} and
Theorem~\ref{th:circ2}.
\end{proof}

Theorem~\ref{th1} and Corollary~\ref{corcirc} provide two different ways to
approximate the solution of biinfinite Toeplitz systems. Which of the two
is preferable? This depends on two criteria:
(i) The accuracy of the approximation for given dimension $n$;
(ii) the computational costs for solving each of the finite-dimensional
systems.

The solution of the circulant system in Corollary~\ref{corcirc} can be
computed via 3 FFTs of size $2n$. The Toeplitz system in Theorem~\ref{th1} 
can be solved by the conjugate gradient method in approximately $3k$ FFTs 
of size $2n$, where $k$ is the number of iterations. Of course, additional
preconditioning can significantly reduce this number at the cost of
two additional FFTs per iteration (see also Section~\ref{s:precond}). For both, 
the circulant and the Toeplitz system, zeropadding can be used to extend 
the vectors to ``power-of-two''-length.

\begin{figure}
\begin{center}
\epsfig{file=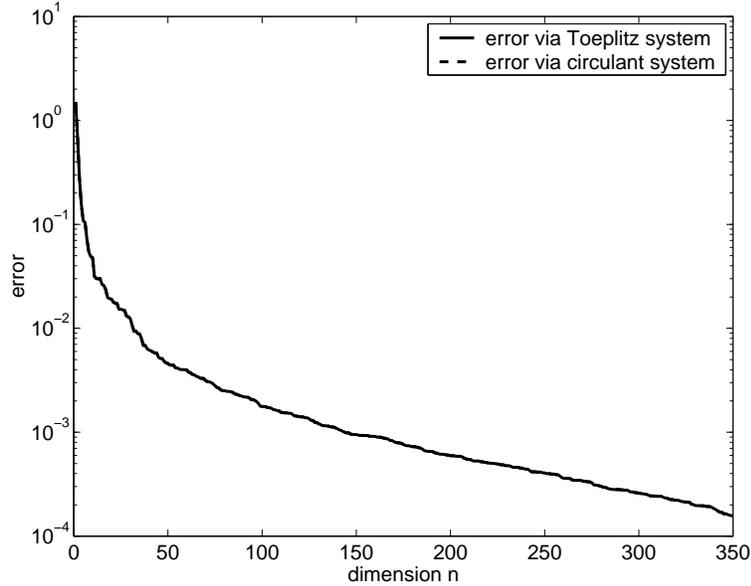,width=100mm}
\caption{Comparison of error for the solution of a biinfinite Toeplitz
system with polynomial decay. We compare the approximation error of the
Toeplitz system described in Theorem~\ref{th1} to that of the circulant 
system of Corollary~\ref{corcirc} for increasing matrix dimension.
The approximation error of both methods is almost identical, so that
the difference between the two graphs is hardly visible.}
\label{fig:fig2}
\end{center}
\end{figure}

{\sc Example 2:}
We consider the same biinfinite Toeplitz system as in Example~1.
We compare the error when approximating the solution by using the circulant 
system of Corollary~\ref{corcirc} and by the Toeplitz system of
Theorem~\ref{th1}. We compute for each $n=1,\dots, 350$ the
approximation error $\|x-x^{(n)}\|$ and $\|x-z^{(n)}\|$ respectively. 
As can be seen from Figure~\ref{fig:fig2} both
methods give almost the same error, in fact the two lines showing the
error are hardly distinguishable. A similar behavior can be observed for other
examples involving biinfinite Toeplitz matrices with fast decay. Since
solving a circulant system is cheaper than solving a Toeplitz system,
the approximation scheme of Corollary~\ref{corcirc} seems to be preferable 
in such situations.

Many variations of the theme are possible. For instance if $A_n$
is an $s$-banded (biinfinite) Toeplitz matrix with $s <n/2$, we could 
use Strang's preconditioner as approximate inverse. Due to the explicit 
constants in Theorem~\ref{th:demko} this approach allows us to give an 
error estimate with explicit constants (cf.~also Theorem~5 in~\cite{Str98a}).
\begin{theorem}
Let $A_n \xn=\yn$ be given where $A_n$ is an $n \times n$ hermitian 
$s$-banded Toeplitz matrix with $s<n/3$ and positive generating function 
and let $\yn_k=0$ for 
$|n/2-k|>s$. Let $C_n$ be the $n\times n$ circulant matrix with first row 
given by $(a_{0},\conj{a}_1,\dots,\conj{a}_s,0,\dots,0, {a}_s,\dots,a_{1})$ 
and let $\zn$ be the solution of $\Cn \zn =\yn$. Then
\begin{equation}
\|\xn-\zn\| \le 3\sqrt{2} c \lam^{-\gamma n} 
(\lam^{-\gamma s}-\lam^{-\gamma(s+1)})^{-3}
\notag 
\end{equation}
where $c$ and $\lam$ are as in Theorem~\ref{th:demko}.
\end{theorem}

\begin{proof}
The proof is similar to that of Theorem~\ref{th:circ2}.
To avoid unnecessary repetitions we only indicate the modifications,
that are required.

By Remark~\ref{remark3} $\Cn$ is invertible. Note that $C_n$ is a matrix 
with three bands, one band is centered at the main diagonal, and the two 
other bands of width $2s$ are located at the lower left and 
upper right corner of the matrix.
It follows from Proposition~5.1 in~\cite{DMS84} that the entries of
$C_n^{-1}$ decay exponentially off the diagonal and off the lower right and
upper left corner. More precisely,
\begin{equation}
|(C_n^{-1})_{k,l}| \le
\begin{cases}
c \lam^{|k-l|} & \text{if $0 \le |k-l| \le n$} \\
c \lam^{2n+1-|k-l|} & \text{if $n+1 \le |k-l| \le 2n$} \,,
\end{cases}
\notag 
\end{equation}
where $c$ is as in Theorem~\ref{th:demko} with $\lam = q^\frac{1}{2s}$.
With this result at hand it is easy to show that the entries of $C_n^{-1} \yn$ 
decay exponentially.

$A_n-C_n$ has a simple form, it is a Toeplitz matrix with
first row 
$$(0,\dots , 0, a_{s} a_{s-1},\dots,a_1).$$

When we compute $\errn:=(A_n-C_n)(C_n^{-1} \yn)$ the non-zero entries of 
$(A_n-C_n)$ are multiplied by exponentially decaying entries
due to the exponential decay of $((C_n^{-1} \yn)$.
This leads to the estimate
$\|\errn\| \le 2 \sqrt{2} c (\lam^{s}-\lam^{s+1})^{-3} \lam^{n}$.
\end{proof}

An interesting alternative to periodic boundary conditions is the use of 
Neumann boundary conditions considered in~\cite{NCT99}. This modification
will be discussed elsewhere.

\section{Preconditioning by embedding and exponentially decaying 
Toeplitz matrices} 
\label{s:precond}

The accuracy of the solution of the Toeplitz system $A_n \xn =\yn$ by using 
$\Sni$ as approximate inverse depends crucially
on the decay properties of the right hand side $\yn$. If $\yn$ does not
have appropriate decay conditions the approach in section~\ref{s:deconv}
may not yield an approximation with sufficient accuracy. But we can still 
use $\Sni$ as preconditioner and solve $\An \xn=\yn$ by the
preconditioned conjugate gradient method~\cite{CN96}. The construction
of preconditioners via circulant embedding is well known, it has
been thoroughly investigated in~\cite{CN93} and for the special case
of band Toeplitz matrices in~\cite{NPT96,HN94}.

Compared to the certainly more elegant and more general approach 
in~\cite{CN93}, the approach undertaken in this section has the advantage 
that it yields some quantitative results. It shows that the clustering 
behavior of the preconditioned matrix $\Sni \An$ is the stronger the faster 
the decay of $\An$ is. Moreover, our approach allows us to prove a conjecture 
by Nagy et al., cf.~\cite{NPT96}, and will provide a theoretical explanation 
for some numerical results presented in \cite{NPT96} and~\cite{HN94}.

The theoretical analysis of the clustering behavior of the eigenvalues
of the preconditioned matrix $\Sni \An$ is inspired by the work of
Raymond Chan~\cite{Cha89a}. We will show that $\Sni \An$ can be written
as $\Sni \An = I_n + R_n +K_n$, where $I_n$ is the identity matrix, $R_n$ 
is a matrix of small rank, and $K_n$ is a matrix of small 2-norm.

\begin{theorem}
\label{th:cluster}
Let $\An$ be a hermitian Toeplitz matrix whose entries $a_k$ decay
exponentially, i.e., $|a_k| \le C e^{\gam |k|}, k=0,\dots,n-1$.
Set $\Sn = \An - \Bn \Ani \Bna$, where $\Bn$ is as defined in~\eqref{embed}.
Then for all $\eps > 0$, there exist $N$ and $M$ such that for all
$n > N$ at most $M$ eigenvalues of $\An - \Sn$ have
absolute value exceeding $\eps$.
\end{theorem}
\begin{proof}
By definition of $\Sn$ we have
\begin{equation}
\An - \Sn = \En,
\notag 
\end{equation}
where $\En= \Bn \Ani \Bna$.
We know from equation~\eqref{bab} in the proof of Theorem~\ref{th:circ2} 
that the entries of $\En$ can be bounded by
\begin{equation}
|(\En)_{kl}| \le 
c (e^{-2\gami n} e^{\gami (k+l)} +  e^{-2\gami n} e^{\gami |k-l|}
+  e^{-\gami(k+l)}).
\label{est3}
\end{equation}

For $x = [x_0,x_1,\dots,x_{n-2},x_{n-1}] \in \Cst^{n}$ and $N<n/2$ we define 
the orthogonal projection $P_N$ by
$$P_N x = [0,\dots,0,x_N,x_{N+1},\dots,x_{n-N-1},0,\dots,0],$$
and identify the image of $P_N$ with $\Cst^{n-2N}$.
We set $\EnN = P_N \En P_N$. In words, 
$\EnN$ is obtained from $\En$ by taking only the central 
$(n-2N) \times (n-2N)$ submatrix of $\En$ and setting the other entries
surrounding this block equal to zero. Then $\En - \EnN$ has $2N$ ``full''
rows and $n-2N$ ``sparse'' rows, where each of the latter rows has non-zero
entries only at the first $N$ and last $N$ coordinates. Thus the dimension
of the space spanned by the sparse rows is at most $2N$. Hence
$\rank (\En - \EnN) \le 4N$.
Due to the decay properties of $\En$ it is easy to see that 
$$\|\EnN\|_1 = \sum_{l=N}^{n-N-1} |(\EnN)_{Nl}|.$$
Using~\eqref{est3} we get after some straightforward calculations 
\begin{equation}
\|\EnN\|_1 \le c (e^{-\gami n}+ e^{-\gami(n-2N)} +e^{-2 \gami N}).
\label{est2}
\end{equation}
It is obvious that for each given $\eps>0$ we can find an $N$ such that for 
all $n>2N$ $\|\EnN\|_1 \le \eps$.
Since $\EnN$ is hermitian, we have $\|\EnN\|_1 = \|\EnN\|_{\infty}$. Thus
$$\|\EnN\|_2 \le (\|\EnN\|_1 \|\EnN\|_{\infty})^{1/2} \le \eps .$$
Hence for large $n$ the spectrum of $\EnN$ lies in $(-\eps,\eps)$. By the 
Cauchy interlace theorem we conclude that at most $4N$ 
eigenvalues of $\An - \Sn$ have absolute value exceeding $\eps$.
\end{proof}

Lemma~\ref{le:cond} implies that for any $\eps >0$ we can find an $M$ such 
that for all $n>M$ $\Sn$ and $\Sni$ exist and $\|\Sni\|$ is bounded by 
$f_{\min}-\eps >0$. Proceeding as in~\cite{Cha89a}, Chapter 2, we express 
$\Sni \An$ as
\begin{equation}
\Sni \An = \In + \Sni (\An - \Sn)
\notag 
\end{equation}
and arrive at the following
\begin{corollary}
\label{cor:chan}
Let $A=\{a_{k,l}\}$ be a hermitian positive definite Toeplitz matrix with
$|a_{k}| \le c e^{-\gamma |k|}$ for $\gamma >0$. Then for all
$\eps >0$ there exist $N$ and $M>0$ such that for all $n>M$ at most
$N$ eigenvalues of $\Sni \An-\In$ have absolute values larger than $\eps$.
\end{corollary}

With the results presented in this paper it should not be difficult
for the reader to derive Theorem~\ref{th:cluster} and
Corollary~\ref{cor:chan} for Toeplitz matrices with polynomial decay.
Theorem~\ref{th:cluster} (in particular~\eqref{est2}) and 
Corollary~\ref{cor:chan} show that the 
clustering behavior of $\Sni \An$ is the stronger the faster the decay 
of the entries of the Toeplitz matrix is.

In~\cite{NPT96} and in~\cite{HN94} Nagy et al.\ consider the
solution of convolution equations and banded Toeplitz systems by
the preconditioned conjugate gradient method using a preconditioner
similar to the one in this section. In~\cite{NPT96} the Toeplitz matrix is
rectangular, but the embedding is in principle the same. The theoretical 
results presented in~\cite{NPT96,HN94} only hold for banded Toeplitz systems, 
but in the numerical experiments the authors consider also
non-banded Toeplitz systems, where the Toeplitz
matrix has exponential decay (see~\cite{NPT96}) or polynomial decay
(see Example~4 in~\cite{HN94}). It is noted in~\cite{NPT96}
that ``it is surprising that the number of iterations is still quite
small''. The authors also point out that the numerical experiments
indicate additional clustering around one of the spectrum of the 
preconditioner matrix, which is not covered by their theoretical
results. With Theorem~\ref{th:cluster} and Corollary~\ref{cor:chan} at hand
we can provide a theoretical explanation for the numerical
observations in~\cite{NPT96} and~\cite{HN94}, at least for the 1-D case.
The key lies in the fast decay of the inverse of the Toeplitz matrix.

\section*{Acknowledgement} 

I want to thank Thomas Kailath and Ali Sayed for fruitful
discussions on this topic, Torsten Ehrhardt for pointing out
reference~\cite{GRS64} to me, Gabriele Steidl for communicating to
me part of the proof of Theorem~\ref{th:circulant} and the referees
for a careful reading of the manuscript.

\end{document}